\documentclass[12pt]{article}

\def\UseBibLatex{1}

\makeatletter
\def\input@path{{styles/}}
\makeatother

\newcommand{\UsePackage}[1]{%
  \IfFileExists{styles/#1.sty}{%
      \usepackage{styles/#1}%
   }{%
      \IfFileExists{../styles/#1.sty}{%
         \usepackage{../styles/#1}%
      }{%
         \usepackage{#1}%
      }%
   }%
}

\usepackage[T1]{fontenc}
\usepackage{lmodern}
\usepackage{textcomp}

\usepackage{amsmath}%
\usepackage{amssymb}%
\usepackage[table]{xcolor}%

\usepackage{todonotes}
\usepackage[in]{fullpage}%

\usepackage[amsmath,thmmarks]{ntheorem}%
\theoremseparator{.}%

\usepackage{titlesec}%
\titlelabel{\thetitle. }%
\usepackage{xcolor}%
\usepackage{mleftright}%
\usepackage{xspace}%
\usepackage{graphicx}
\usepackage{hyperref}%

\usepackage{hyperref}%
\hypersetup{%
      unicode,
      breaklinks,%
      colorlinks=true,%
      urlcolor=[rgb]{0.25,0.0,0.0},%
      linkcolor=[rgb]{0.5,0.0,0.0},%
      citecolor=[rgb]{0,0.2,0.445},%
      filecolor=[rgb]{0,0,0.4},
      anchorcolor=[rgb]={0.0,0.1,0.2}%
}
\usepackage[ocgcolorlinks]{ocgx2}

\theoremseparator{.}%

\theoremstyle{plain}%
\newtheorem{theorem}{Theorem}[section]

\newtheorem{lemma}{Lemma}[section]

\newtheorem{definition}{Definition}[section]

\theoremstyle{plain}
\theoremheaderfont{\bfseries}
\theorembodyfont{\upshape}
\newtheorem{remark}{Remark}[section]

\theoremstyle{plain}%

\newcommand{\thmlab}[1]{\label{theo:#1}}
\newcommand{\thmref}[1]{\HLink{Theorem}{theo:#1}}

\newcommand{\lemlab}[1]{\label{lemma:#1}}
\newcommand{\lemref}[1]{\HLink{Lemma}{lemma:#1}}

\providecommand{\deflab}[1]{}
\renewcommand{\deflab}[1]{\label{def:#1}}
\newcommand{\defref}[1]{\HLink{Definition}{def:#1}}

\newcommand{\remlab}[1]{\label{rem:#1}}
\newcommand{\remref}[1]{\HLink{Remark}{rem:#1}}

\providecommand{\emphind}[1]{}%
\renewcommand{\emphind}[1]{\emph{#1}\index{#1}}

\definecolor{blue25emph}{rgb}{0, 0, 11}

\providecommand{\emphic}[2]{}
\renewcommand{\emphic}[2]{\textcolor{blue25emph}{%
      \textbf{\emph{#1}}}\index{#2}}

\providecommand{\emphi}[1]{}%
\renewcommand{\emphi}[1]{\emphic{#1}{#1}}

\definecolor{almostblack}{rgb}{0, 0, 0.3}

\providecommand{\emphw}[1]{}%
\renewcommand{\emphw}[1]{{\textcolor{almostblack}{\emph{#1}}}}%

\providecommand{\emphOnly}[1]{}%
\renewcommand{\emphOnly}[1]{\emph{\textcolor{blue25}{\textbf{#1}}}}

\newcommand{\HaoranThanks}[1]{%
   \thanks{%
       Department of Mathematics; %
      Northeastern University; %
      567 Lake Hall;%
      43 Leon Street; %
      Boston, MA, 02115, USA; %
      \href{wu.haoran2@northeastern.edu}{wu.haoran2@northeastern.edu}; %
      \url{}. %
   #1%
   }%
}

\newcommand{\HLink}[2]{\hyperref[#2]{#1~\ref*{#2}}}
\newcommand{\HLinkSuffix}[3]{\hyperref[#2]{#1\ref*{#2}{#3}}}

\providecommand{\deflab}[1]{}
\renewcommand{\deflab}[1]{\label{def:#1}}

\renewcommand{\qed}{\hfill$\square$}
\providecommand{\eqlab}[1]{}%
\renewcommand{\eqlab}[1]{\label{equation:#1}}

\newcommand{\remove}[1]{}%

\usepackage[inline]{enumitem}

\newlist{compactenumA}{enumerate}{5}%
\setlist[compactenumA]{topsep=0pt,itemsep=-1ex,partopsep=1ex,parsep=1ex,%
   label=(\Alph*)}%

\newlist{compactenuma}{enumerate}{5}%
\setlist[compactenuma]{topsep=0pt,itemsep=-1ex,partopsep=1ex,parsep=1ex,%
   label=(\alph*)}%

\newlist{compactenumI}{enumerate}{5}%
\setlist[compactenumI]{topsep=0pt,itemsep=-1ex,partopsep=1ex,parsep=1ex,%
   label=(\Roman*)}%

\newlist{compactenumi}{enumerate}{5}%
\setlist[compactenumi]{topsep=0pt,itemsep=-1ex,partopsep=1ex,parsep=1ex,%
   label=(\roman*)}%

\newlist{compactitem}{itemize}{5}%
\setlist[compactitem]{topsep=0pt,itemsep=-1ex,partopsep=1ex,parsep=1ex,%
   label=\ensuremath{\bullet}}%

\providecommand{\BibLatexMode}[1]{}
\providecommand{\BibTexMode}[1]{}

\ifx\UseBibLatex\undefined%
  \renewcommand{\BibLatexMode}[1]{}
  \renewcommand{\BibTexMode}[1]{#1}
\else
  \renewcommand{\BibLatexMode}[1]{#1}
  \renewcommand{\BibTexMode}[1]{}
\fi

\BibLatexMode{%
   \usepackage[bibencoding=utf8,style=numeric,citestyle=numeric,backend=biber,sorting=none,]{biblatex}%
   \UsePackage{my_biblatex}%
}

\numberwithin{figure}{section}%
\numberwithin{table}{section}%
\numberwithin{equation}{section}%

\usepackage{graphicx} 

\usepackage{chngcntr}
\counterwithout{figure}{section}

\newcommand{\appref}[1]{\hyperref[#1]{\textcolor{red}{Appendix}}}

\DeclareUnicodeCharacter{2212}{-}
\date{}  
\begin{document}

\title{Euler equation on a fast rotating ellipsoid }

\author{Haoran Wu\HaoranThanks{}}

\maketitle

\begin{abstract}

This paper extends the analytical study of the incompressible Euler equations from the classical spherical setting to the more realistic geometry of a biaxial ellipsoid. Motivated by the work of Cheng and Mahalov \cite{Cheng2013} on fast rotating spheres and Xu \cite{Xu2024} on Rossby–Haurwitz solutions on ellipsoids, we adapt their framework to establish a parallel result for Euler flows on a rotating ellipsoidal surface. In the regime of rapid rotation, we prove that the time–averaged velocity field remains uniformly bounded in Sobolev norms independent of the rotation rate and converges to a longitude–independent zonal flow depending only on latitude. This shows that the zonalization phenomenon discovered by Cheng and Mahalov on the sphere persists on biaxial ellipsoids, thereby bridging the gap between spherical and ellipsoidal theories of fast rotating Euler dynamics.

\end{abstract}


\section{Introduction}



One of central topics in geophysical fluid dynamics is the dynamics of incompressible Euler flows on rotating planetary surfaces. The spherical geometry has served as the canonical setting, both for analytical study and for numerical simulation. In particular, Cheng and Mahalov\cite{Cheng2013} established that finite-time averages of Euler flows on a fast rotating sphere tend toward zonal states. 

However, real planets are not perfect spheres. Most of them such as Jupiter or Saturn introduces a natural deviation toward ellipsoidal geometry, where the variation in curvature and the modified Coriolis profile can significantly influence the fluid dynamics. A fast rotating Jovian planet usually has a relatively large flattening rate (see Berardo and Wit \cite{BerardoWit2022}). Saturn has a large flattening rate about $0.1$ (see Elkins–Tanton \cite{ElkinsTanton2006}), and Haziot \cite{Haziot2022} have shown that a spherical model turned to be unsuitable for flows on Saturn. Hence, a biaxial ellipsoid model is required to more faithfully capture the geometric characteristics of outer planets. Recent studies have begun to address this more realistic geometry. For instance,
Taylor \cite{Taylor2016} investigated the two-dimensional Euler equations on a general rotationally symmetric manifold to better understand the dynamics of planetary atmospheres, with particular attention to the stability of zonal flow solutions. Inspired by these works, this article investigates the 2D Euler equation on a fast rotating ellopsoid. Our approach is to adapt the PDE–geometric framework of Cheng–Mahalov\cite{Cheng2013}, originally formulated for the fast rotating sphere, to the setting of a biaxial ellipsoid built by \cite{Xu2024}. By following the same decomposition methods, time-averaging techniques, and operator estimates, we investigate how the geometry of the ellipsoid modifies the structure of the Euler equations and the associated invariant subspaces. In particular, we show how the Hodge decomposition, the identification of the null space associated with the Coriolis operator, and the key Sobolev estimates can be reformulated in the ellipsoidal geometry. This yields an analytical description of the leading-order dynamics and clarifies the role of the meridional dependence of the Coriolis parameter in the absence of full spherical symmetry.

In summary, this paper extends the analytic framework of Cheng and Mahalov for time-averaged Euler flows on the sphere to the biaxial ellipsoid. We establish uniform Sobolev estimates independent of the rotation rate and prove that the time-averaged flow converges to a latitude-dependent zonal state, thereby clarifying how geometric flattening modifies the structure of fast-rotating incompressible flows. Our contribution provides a step toward a more accurate mathematical understanding of atmospheric flows on rotating planets with significant flattening, aligning theoretical analysis with physical observations

In this study, we work on a biaxial ellipsoid $\mathbb{E}^2$ with major axis equal to 1 and minor axis $b<1$, rotating with angular velocity $\omega$, see Fig 1\cite{Xu2024}. Following the same setting in \cite{Xu2024}, a point $P$ on its surface is parameterized by the angular coordinates $(\varphi,\theta)\in[-\pi,\pi)\times\left[-\frac{\pi}{2},\frac{\pi}{2}\right]$, where $(x, y, z) = (\cos \theta \cos \varphi, \cos \theta \sin \varphi, b \sin \theta)$. The equator corresponds to $\theta = 0$, the North Pole to $\theta = \frac{\pi}{2}$, and the South Pole to $\theta = -\frac{\pi}{2}$.

\begin{figure}[htbp]
  \centering
\includegraphics[width=1\textwidth]{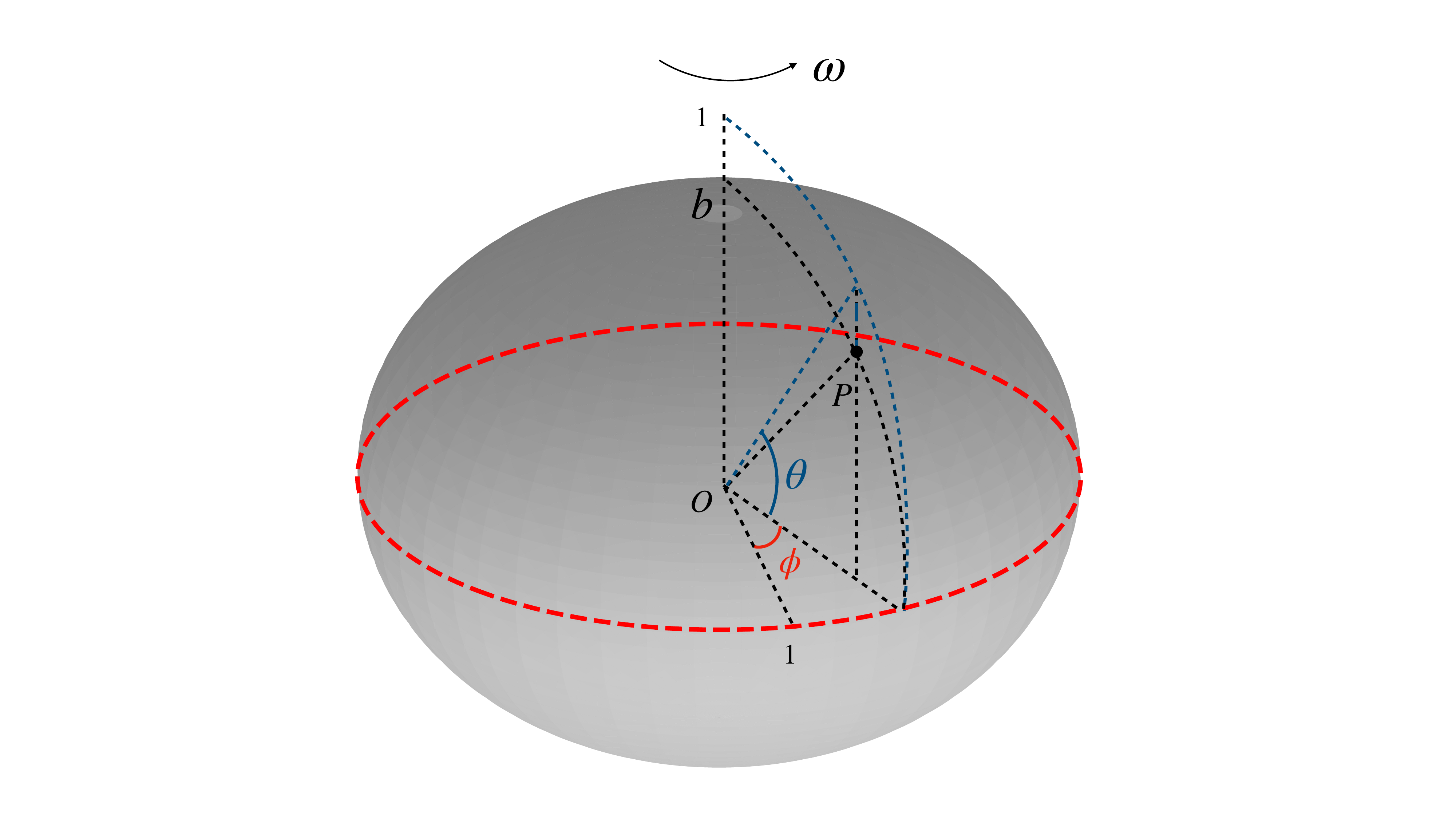}  
\caption{A biaxial ellipsoid rotating with angular velocity $\omega$\cite{Xu2024}}
  \label{fig:Ellipsoid}
\end{figure}

The Euler equations under such coordinate in terms of the velocity field $u$ is 

\begin{equation}
    \begin{cases}
D_{t} u 
+ \dfrac{2\omega \sin\theta}{\sqrt{\sin^{2}\theta + b^{2}\cos^{2}\theta}} \, J u 
= - \nabla p, \\[1.2ex]
\operatorname{div} u = 0.
\end{cases} \label{1.1}
\end{equation}
where $D_{t} u=\partial_t u + \nabla_{u}u$, $\omega$ is the rotating speed of the ellipsoid, $J$ is a rotation matrix, i.e., $Ju=u^{\perp}$ where $u^{\perp}$ denotes a counter-clockwise $\pi/2$-rotation of $u$ on $\mathbb{E}^2$. And  $p$ is the pressure field determined implicitly by the velocity field $u$. By simplifyinging, we have the following.

\begin{equation}
    \partial_tu+\nabla p=\frac{-2w\sin\theta}{\sqrt{\sin^2\theta+b^2\cos^2\theta}}u^{\perp}-\nabla_{u}u \label{eulerequation}
\end{equation} 

\noindent Our theoretical investigation is then focused on the fast rotating
regime with $\varepsilon\ll1$. $\varepsilon$, called the Rossby number, defines as $\varepsilon = U/2\omega \sin\theta$, $U$ and $L$ respectively characteristic velocity and length scales of the flow. It scales like the frequency of the frame’s rotation and when $\varepsilon\ll1$, the Coriolis force dominates and the flow tends to be zonal. Under such setting, the nature of the time-averages of $u$:

\begin{equation}
    \overline{{u}}(T,\cdot):=\frac1T\int_0^T {u}(t,\cdot)dt \label{1.2}  
\end{equation}
The main result is stated as the following.

\begin{theorem}\thmlab{thm1.1}

Consider the incompressible Euler equation \eqref{1.1} on $\mathbb{E}^2$ with initial data $u_0\in H^k(\mathbb{E}^2)$ for $k\geq3$. Define the time-averaged flow $\overline{{u}}$ as in \eqref{1.2}. Then, there exists a function $f(\cdot):[-1,1]\mapsto\mathbb{R}$ and constants $C_0,T_0$ independent of $\omega$ and $u_0$,s.t.for any given $T\in[0,T_0/\|{u}_0\|_{H^k}]$,

\begin{equation}
    \left\|\overline{{u}}(T,x,y,z)-\nabla^\perp f(z)\right\|_{H^{k-3}(\mathbb{E}^2)}\leq \frac{C_0}{\omega}\left(\frac{M_0}T+M_0^2\right) \label{1.3}
\end{equation}
with $M_0:=\|{u}_0\|_{H^k}$.In elliptic coordinates, the approximation $\nabla^\perp f(z)$ is

$$\nabla^\perp f(z)=-f^{\prime}(\frac{-2\sin\theta}{\sqrt{\sin^{2}\theta+b^{2}\cos^{2}\theta}})\frac{-2b^{2}\cos\theta}{(\sin^{2}\theta+b^{2}\cos^{2}\theta)^{2}}e_{\phi}$$
which is a longitude-independent zonal flow.
   
\end{theorem}
The proof of our main theorem starts from \eqref{eulerequation} using Hodge decomposition into 
\begin{equation} \partial_{t}u+\nabla^{\perp}\Delta^{-1}\mathrm{curl}(\nabla_{u}u)=w\cdot \mathcal{L}[u]
\end{equation}
where $\mathcal{L}[u]=\nabla^{\perp}\Delta^{-1}\mathrm{curl}(\frac{-2\sin\theta}{\sqrt{\sin^{2}\theta+b^{2}\cos^{2}\theta}}u^{\perp})$. Then, the large constant $\omega$ leads to an $O(\omega)$ estimate on the time-average of the $\mathcal{L}[u]$ term and eventually leads to the Main Theorem. This procedure fits into the abstract framework of the following \lemref{lemma1.1}. 

Our \lemref{lemma1.1} is parallel to the corresponding averaging lemma in Cheng and Mahalov\cite{Cheng2013}.
The structure and proof are essentially identical, except that we replace $\varepsilon$ with our notation $\omega$. For completeness, we restate the lemma below in our current setting on the ellipsoidal surface. Detailed proofs and intermediate estimates can be found in Cheng and Mahalov\cite{Cheng2013}[Lemma 1.1], and the same argument applies verbatim here with $1/\omega$ replacing $\varepsilon$.

\begin{lemma}\cite{Cheng2013}\lemlab{lemma1.1} Consider time-dependent equation, 

$$\partial_t{u}=\omega\mathcal{L}[{u}]+f$$
where $\varepsilon=U/2\omega L$, and $0<\varepsilon\ll1$, $\mathcal{L}$ is a linear operator and $f$ includes nonlinear and source terms. Assume a priori that, for some Hilbert spaces $X_1$,$X_2$, one has

\begin{equation}
{u}\in\mathcal{C}([0,T],{X}_1\cap{X}_2),\quad\mathcal{L}[{u}]\in\mathcal{C}([0,T],{X}_2),\quad f\in\mathcal{C}([0,T],{X}_2) \label{1.5}
\end{equation}
and for such solution, 

\begin{equation}
\int_0^T\mathcal{L}[{u}]dt=\mathcal{L}\left[\int_0^T{u}dt\right] \label{1.6}
\end{equation}
Also, let operator $\prod _{{null}\{ \mathcal{L} \} }: {X} _1\to {X} _1$ denote (some) projection onto the kernel of $\mathcal{L}$. Then, under the assumption

\begin{equation}
    \|{u}-\prod_{\mathrm{null}\{\mathcal{L}\}}{u}\|_{{X}_1}\leq C\|\mathcal{L}[{u}]\|_{{X}_2} \label{1.7}
\end{equation}
for some constant $C$, the following estimate holds true on the time—average of $u$,

$$\left\|\frac1T\int_0^T{u}\:dt-\frac1T\int_0^T\prod_{\mathrm{null}\{\mathcal{L}\}}{u}\:dt\right\|_{{X}_1}\leq\frac{C}{\omega}\left(\frac{2M}T+M^{\prime}\right)$$
where constants $M:=\max_{t\in[0,T]}\|{u}(t,\cdot)\|_{{X}_2}$ and $M^\prime:=\max_{t\in[0,T]}\|f(t,\cdot)\|_{X_2}$   
\end{lemma}

\begin{remark}
    The idea of estimating time-averages for PDE systems with fast oscillations has appeared in e.g. \cite[Th.~2.5]{Cheng2012_SIAM},~\cite{ChengMahalov2012_DCDS}.
\end{remark}
We will use this lemma to organize our proof for  \thmref{thm1.1}. In a nutshell, we will define operators $\mathcal{L}[{u}]:= \nabla^{\perp}\Delta^{-1}\mathrm{curl}(\frac{-2\sin\theta}{\sqrt{\sin^{2}\theta+b^{2}\cos^{2}\theta}}u^{\perp})$ in \defref{def2.1} and $\prod_{\mathrm{null}\{\mathcal{L}\}}$ in \lemref{lem3.2} prove the key estimate \eqref{1.7} in \thmref{thm4.1} and verify the assumptions \eqref{1.5} and \eqref{1.6} in the last section.

\section{Hodge decomposition}
The Hodge decomposition theorem \cite{Taylor1996,Warner1983} confirms that for any $k$-form $\omega$ on an oriented compact Riemannian manifold, there exist a $(k-1)$-form $\alpha$, $(k+1)$-form $\beta$ and a harmonic $k$-form $\gamma$, s.t. $\omega=d\alpha+\delta\beta+\gamma$. In particular, if the manifold is the surface in the cohomology class, then for any smooth vector field $u$ we could decomposite it into two part, represented by two scalar-valued functions $\Phi$ (called potential) and $\Psi$ (called stream function), such that

\begin{equation}
    {u}={u}_{\mathrm{irr}}+{u}_{\mathrm{inc}}\quad\text{where~}{u}_{\mathrm{irr}}:=\nabla\Phi\mathrm{~and~}{u}_{\mathrm{inc}}:=\nabla^\perp\Psi \label{hodgeproperty1}
\end{equation}

\noindent We use subscript ‘irr’ stands for irrotational vector fields and ‘inc’ for incompressible. In other word, any such manifold could be divided into irrotational part and incompressible part uniquely . Moreover, we have such property

\begin{equation}
    \mathrm{curl}~{u}_{\mathrm{irr}}=\mathrm{div}~{u}_{\mathrm{inc}}=0,\quad{u}_{\mathrm{irr}}=\nabla\Delta_{\mathbb{E}^2}^{-1}\mathrm{div}~{u} \quad\text{and}\quad{u}_{\mathrm{inc}}=\nabla^\perp\Delta_{\mathbb{E}^2}^{-1}\mathrm{curl}~{u} \label{hodgeproperty}
\end{equation}

\noindent Observe \eqref{eulerequation} on the RHS, $\partial_tu$ is incompressible and $\nabla p$ is irrotational. Thus, the RHS is the unique Hodge decomposition of the
LHS. Then we could use the property above \eqref{hodgeproperty}. In particular,the incompressible part $\partial_tu$ is uniquely determined by

\begin{equation}
    \partial_{t}u=\nabla^{\perp}\Delta^{-1}\mathrm{curl}~(\frac{-2w\sin\theta}{\sqrt{\sin^{2}\theta+b^{2}\cos^{2}\theta}}u^{\perp}-\nabla_{u}u) \label{partialtu}
\end{equation}

\noindent In the context of \lemref{lemma1.1}, we define the following operator.
\begin{definition}\deflab{def2.1}

For any vector field $u$ on $\mathbb{E}^2$, not necessarily div-free, we define

\begin{equation}
    \mathcal{L}[u]=\nabla^{\perp}\Delta^{-1}\mathrm{curl}~(\frac{-2\sin\theta}{\sqrt{\sin^{2}\theta+b^{2}\cos^{2}\theta}}u^{\perp}) \label{Loperator}
\end{equation}
thus \eqref{partialtu} can be reformulated as 

\begin{equation}
    \partial_{t}u+\nabla^{\perp}\Delta^{-1}\mathrm{curl}~(\nabla_{u}u)=w\cdot \mathcal{L}[u] \label{partialtubyLopertaor}
\end{equation}
\end{definition}

\section{Null space of L and associated $L^2$-orthogonal projection}

By observing \eqref{Loperator}, it is easy to obtain that the sufficient condition for $\mathcal{L}[u]=0$ is $$\mathrm{curl}~(\frac{-2\sin\theta}{\sqrt{\sin^{2}\theta+b^{2}\cos^{2}\theta}}u^{\perp})=0$$ By property of $\mathrm{curl}$, i.e., $\mathrm{curl}~\nabla^\perp=\Delta$ it is also easy to verify the necessity. Thus we have 

\begin{equation}
    \mathcal{L}[{u}]=0\Longleftrightarrow\mathrm{~curl~}(\frac{-2\sin\theta}{\sqrt{\sin^{2}\theta+b^{2}\cos^{2}\theta}}u^{\perp})=0
\end{equation}

\begin{lemma}\lemlab{lemma3.1} (Characterization of null $\mathcal{L}$) For div-free vector field $u$ on $\mathbb{E}^2$ with sufficient regularity,

\begin{equation}
    \mathcal{L}[u]=0\Leftrightarrow
 u=\Psi^{\prime}(\frac{-2\sin\theta}{\sqrt{\sin^{2}\theta+b^{2}\cos^{2}\theta}})\frac{-2b^{2}\cos\theta}{(\sin^{2}\theta+b^{2}\cos^{2}\theta)^{2}}e_{\phi}
\end{equation}
for some function $\Psi :[-1,1]\mapsto\mathbb{R}$. Thus we identify  $\mathrm{null}\{\mathcal{L}\}$ when restricted to div-free velocity fields, with the space of longitude−independent zonal flows.
\end{lemma}

\noindent\textbf{Proof of Lemma 3.1.} By property of $\mathrm{curl}$, $\mathrm{curl}~\left(f{u}^\perp\right)=\mathrm{div}~\left(f{u}\right)=\nabla f\cdot{u}+f\mathrm{div}~{u}$. we have:

$$\mathrm{curl}~(\frac{-2\sin\theta}{\sqrt{\sin^2\theta+b^2\cos^2\theta}}u^{\perp})=\left(\nabla\frac{-2\sin\theta}{\sqrt{\sin^{2}\theta+b^{2}\cos^{2}\theta}}\right)u$$\\
So $\mathcal{L}[u]=0$ implies $(\nabla\frac{-2\sin\theta}{\sqrt{\sin^{2}\theta+b^{2}\cos^{2}\theta}})\cdot u=0$. By Hodge decomposition properties \eqref{hodgeproperty1}, \eqref{hodgeproperty}, we could obtain that $\mathrm{div~}u=0$ is equivalent to $u=\nabla^{\perp}\Psi$. Then we have, for any incompressible velocity field $u$, such that $u=\nabla^{\perp}\Psi \in \mathrm{null}\{\mathcal{L}\}$, we have

$$(\nabla\frac{-2\sin\theta}{\sqrt{\sin^{2}\theta+b^{2}\cos^{2}\theta}})\|(\nabla\Psi)$$
Above condition implies that $\Psi$ is a function of $\theta$ only. By $\mathrm{grad~}\psi=\frac{\partial_\varphi\psi}{\cos\theta}\mathbf{e}_\varphi+\frac{\partial_\theta\psi}{\sqrt{\sin^2\theta+b^2\cos^2\theta}}\mathbf{e}_\theta$, we could have:

$$\nabla\frac{-2\sin\theta}{\sqrt{\sin^2\theta+b^2\cos^2\theta}}=\frac{-2b^2\cos\theta}{(\sin^2\theta+b^2\cos^2\theta)^2}e_\theta$$
and
$$\nabla^{\perp}\frac{-2\sin\theta}{\sqrt{\sin^2\theta+b^2\cos^2\theta}}=\frac{-2b^2\cos\theta}{(\sin^2\theta+b^2\cos^2\theta)^2}e_\phi$$
So $u$ can be reformulated as 

$$u=\nabla^{\perp}\Psi(\frac{-2\sin\theta}{\sqrt{\sin^2\theta+b\cos^2\theta}})=\Psi^{\prime}(\frac{-2\sin\theta}{\sqrt{\sin^2\theta+b\cos^2\theta}})\cdot\frac{-2b^2\cos\theta}{(\sin^2\theta+b^2\cos^2\theta)^2}e_\phi$$

\qed

\noindent Now we need to consider characterization of $\prod_{\mathrm{null}\{\mathcal{L}\}}$, it is easy to show that it is the $L^2$-orthogonal-projection operator onto $\mathrm{null}\{\mathcal{L}\}$.

\begin{lemma}\lemlab{lem3.2}

For any div-free vector field ${u}\in L^2(\mathbb{E}^2)$, it's $L^2$-orthogonal-projection onto ${\mathrm{null}\{\mathcal{L}\}}$ satisfies:

\begin{equation}
    \prod_{\mathrm{null}\{\mathcal{L}\}}u=\prod_{\mathrm{null}\{\mathcal{L}\}}\frac{(\oint_{C(\theta)}u\cdot e_\phi)e_\phi}{\oint_{C(\theta)}e_\phi\cdot e_\phi}=\frac{\sqrt{\sin^{2}\theta+\frac{\cos^{2}\theta}{b^{2}}}}{2\pi\cos\theta}(\oint_{C(\theta)}u\cdot e\phi)e_{\phi}
\end{equation}
where $\oint_{C(\theta)}$ is the line integral along the ellipse at a fixed latitude $\theta$.

\end{lemma}

\noindent\textbf{Proof of Lemma 3.2.} For any $u\in L^2(\mathbb{E}^2)$, it's $L^2$-orthogonal-projection onto ${\mathrm{null}\{\mathcal{L}\}}$ satisfies:
$$\prod_{\mathrm{null}\{\mathcal{L}\}}u=\prod_{\mathrm{null}\{\mathcal{L}\}}\frac{(\oint_{C(\theta)}u\cdot e_\phi)e_\phi}{\oint_{C(\theta)}e_\phi\cdot e_\phi}$$
$\oint_{C(\theta)}e_\phi\cdot e_\phi$ is the length of the ellipse at latitude $\theta$. By the coordinates we set $(x,y,z)=(r(\theta)\sin\theta\cos\varphi,r(\theta)\sin\theta\sin\varphi,r(\theta)\cos\theta)$ and $x^{2}+y^{2}+\frac{z^{2}}{b^{2}}=1$. It is easy to get

$$(\oint_{C(\theta)}u\cdot e_\phi)e_\phi =\frac{2\pi\cdot\cos\theta}{\sqrt{\sin^{2}\theta+\frac{\cos^{2}\theta}{b^{2}}}}$$
\qed

\noindent In terms of stream function, as what we have proved beofore for any div-free velocity field $u=\nabla^{\perp}\Psi$, with amount to $u=(\partial_\theta \Psi)e_{\theta}^{\perp}+(\frac{\partial_\phi\Psi}{\sin\theta})e_{\phi}^{\perp}$. With
removable singularity at the poles,

\begin{equation}
    \begin{aligned}\prod_{\mathrm{null}\{\mathcal{L}\}}(\nabla^{\perp}\Psi)&=\frac{\sqrt{\sin^{2}\theta+\frac{\cos^{2}\theta}{b^{2}}}}{2\pi\cos\theta}\left(\int_{0}^{2\pi}\partial_\theta\Psi d\phi\right)e\phi\\&=\nabla^{\perp}\left(\frac{\sqrt{\sin\theta+\frac{\cos^{2}\theta}{b^{2}}}}{2\pi}\int_{0}^{2\pi}\Psi d\phi\right)\end{aligned} \label{projection}
\end{equation}

\noindent In other words, the operator $\prod_{\mathrm{null}\{\mathcal{L}\}}$ maps the stream function to its zonal means.

\section{Key estimates}
This section is to provide an estimates with the tools of Sobolev norms and inner product associated with the Hilbert space $L^2(\mathbb{E}^2)$. The following lemma and theorem are inspired by the work of Cheng and Mahalov\cite{Cheng2013}, with minor modifications adapted to our ellipsoidal setting. In this context, a natural tool in studying Sobolev norms and $L^2$ inner products on $\mathbb{E}^2$ is the  ellipsoidal harmonics. Not like the sphere case, we know the exact eigenvalue, here we find that our etimates do not reply that much on the exact value, so by our construction we denote the "eigenvalues" for the ellipsoidal harmonics as $\Lambda_{l,m}$, which satisfy $-\Delta Y_l^m=\Lambda_{l,m}Y_l^m$. And we shall still have $\partial_\phi Y_l^m=imY_l^m$ due to the the symmetry that
is preserved in this ellipsoidal case under our consideration. Some spectral theory results (Weyl’s law) ensure that: $|\Lambda_{l,m}|\to \infty$ as $l \to \infty$. Note that the values of $l$ and $m$ for the ellipsoidal harmonics are different from the spherical
harmonics, in fact $l$ are probably no longer integers (although in the perturbative
regime $b\approx1$ they should be close to the spherical case. The ”eigenvalues” $\Lambda_{l,m}$ for the ellipsoidal harmonics are also different but should be close to $l(l+1)$ when  $b\approx1$. We will just use the notation $\sum_{l=1}^\infty\sum_{m=-l}^l$ below for convenience. We set 
\begin{equation}
    \frac{b^2}{(\sin^2\theta+b^2\cos^2\theta)^2}\Psi=\sum_{l,m}\tilde{\psi}_l^mY_l^m \quad\text{with}\quad \tilde{\psi}_l^m=\langle\frac{b^2}{(\sin^2\theta+b^2\cos^2\theta)^2}\Psi,Y_l^m\rangle_{L^2(\mathbb{S}^2)} \label{harmonicforpsi}
\end{equation}

\begin{remark}\remlab{rmk 4.1}
    Here and below, we assume $\psi_0^0=\int_{\mathbb{S}^2}\Psi=0$ and exclude $l=0$ from any series.
\end{remark}

\noindent Following the motivated definition of the $H^k$ space in the Cheng-Mahalov paper\cite{Cheng2013},
define (here A$\thicksim$B means A/B is bounded uniformly from above and below by positive constants that only depend on $k$)

$$\|f\|_{H^k(\mathbb{E}^2)}\thicksim\sqrt{\sum_{j=0}^k(-1)^j\langle\Delta^jf,f\rangle_{L^2(\mathbb{S}^2)}}$$

\noindent And hence, 

$$\|Y_l^m\|_{H^k(\mathbb{E}^2)}\thicksim\sqrt{\sum_{j=0}^k\Lambda_{l,m}^j}\thicksim\sqrt{\Lambda_{l,m}^k}$$
Using the above relations, we are now ready to introduce the following definition.
\begin{definition}
    For a scalar function $\Psi=\sum_{l,m}\psi_l^mY_l^m$ on $\mathbb{E}^2$ with $\int_{\mathbb{E}^2}\Psi=0$, we define its $H^k$ norm as:
$$\|\Psi\|_{H^k(\mathbb{E}^2)}:=\sqrt{\sum_{m,l}\Lambda_{l,m}^k|\psi_l^m|^2}$$
\end{definition}
We can then define $H^k$ norms for a vector field $u$ on $\mathbb{E}^2$.

\begin{definition}
    For a vector field u with Hodge Decomposition ${u}=\nabla\Phi+\nabla^{\perp}\Psi,\int_{\mathbb{E}^{2}}\Phi=\int_{\mathbb{E}^{2}}\Psi=0$, one can define

\begin{equation}
    \|{u}\|_{H^k(\mathbb{E}^2)}:=\sqrt{\|\Phi\|_{H^{k+1}(\mathbb{E}^2)}^2+\|\Psi\|_{H^{k+1}(\mathbb{E}^2)}^2}
\end{equation}
In particular, if $u$ is div-free with ${u}=\nabla^\perp\Psi$ and $\int_{\mathbb{E}^2}\Psi=0$, then

\begin{equation}
    \|\nabla^\perp\Psi\|_{H^k(\mathbb{E}^2)}=\sqrt{\sum_{m,l}\Lambda_{l,m}^{k+1}|\psi_l^m|^2} \label{4.7}
\end{equation}

\end{definition}

\begin{remark}
Here, we follow the zero mean convention in \remref{rmk 4.1}, so that the above definition is consistent with $\|\mathbf{0}\|_{H^k}=0$.
\end{remark}

\noindent With the help of ellipsoidal harmonics we can get the representation for $\mathcal{L}$ operator. From \eqref{harmonicforpsi}, we obtain

$$\frac{b^2}{(\sin^2\theta+b^2\cos^2\theta)^2}\partial_\phi\Psi=\partial_\phi[\frac{b^2}{(\sin^2\theta+b^2\cos^2\theta)^2}\Psi]=\sum_{l=1}^\infty\sum_{m=-l}^lim\tilde{\psi}_l^mY_l^m$$

\noindent Notice that 

\begin{equation}
    \mathcal{L}(\nabla^\perp\Psi)=\nabla^\perp\Delta^{-1}(\frac{b^2}{(\sin^2\theta+b^2\cos^2\theta)^2}\partial_\phi\Psi)=\nabla^\perp\Delta^{-1}\partial_\phi[\frac{b^2}{(\sin^2\theta+b^2\cos^2\theta)^2}\Psi] \label{identity}
\end{equation}

\begin{lemma}\lemlab{lem4.1} (spherical-harmonic representation of $\mathcal{L}$) For a scalar function $\Psi$ defined as \eqref{harmonicforpsi}
combined with identity \eqref{identity}, we can obtain

\begin{equation}
    \mathcal{L}(\nabla^\perp\Psi)=\nabla^\perp\sum_{l=1}^\infty\sum_{m=-l,}^l\frac{-im}{\Lambda(l)}\tilde{\psi}_l^mY_l^m \label{sphericalharmonicrepresentationofL}
\end{equation}

\end{lemma}

\noindent Now we use ellipsoidal harmonics to get representation of the projection operator $\prod_{\mathrm{null}\{\mathcal{L}\}}$ defined in \eqref{projection}, it follows from \eqref{sphericalharmonicrepresentationofL} that if $u=\nabla^{\perp}\Psi\in\mathrm{null}(\mathcal{L})$ then we have 
$$\nabla^{\perp}\sum_{l=1}^{\infty}\sum_{m=-l,}^{l}\frac{-im}{\Lambda(l)}\tilde{\psi}_{l}^{m}Y_{l}^{m}=0$$
by \eqref{harmonicforpsi}, we have 
$$\frac{b^{2}}{(\sin^{2}\theta+b^{2}\cos^{2}\theta)^{2}}\Psi=\sum_{l=1}^{\infty}\tilde{\psi}_{l}^{0}Y_{l}^{0}$$
After simplifying we have the following result
$$\Psi=\frac{(\sin^{2}\theta+b^{2}\cos^{2}\theta)^{2}}{b^{2}}\sum_{l=1}^{\infty}\tilde{\psi}_{l}^{0}Y_{l}^{0}=\sum_{l=1}^{\infty}\frac{(\sin^{2}\theta+b^{2}\cos^{2}\theta)^{2}}{b^{2}}\tilde{\psi}_{l}^{0}Y_{l}^{0}$$
which is consistent with \lemref{lemma3.1}, since $Y_{l}^{0}$ is a function of $\theta$ only.

\begin{lemma} (spherical-harmonic representation of $\prod_{\mathrm{null}\{\mathcal{L}\}}$) For a scalar function $\Psi$ defined as \eqref{harmonicforpsi} we have 

\begin{equation}
    \prod_{\mathrm{null}(\mathcal{L})}(\nabla^\perp\Psi)=\nabla^\perp\sum_{l=1}^\infty\frac{(\sin^2\theta+b^2\cos^2\theta)^2}{b^2}\tilde{\psi}_l^0Y_l^0=\nabla^\perp\frac{(\sin^2\theta+b^2\cos^2\theta)^2}{b^2}\sum_{l=1}^\infty\tilde{\psi}_l^0Y_l^0
\end{equation}

\begin{equation}
    (\mathrm{id}-\prod_{\mathrm{null}(\mathcal{L})})(\nabla^\perp\Psi)=\nabla^\perp\frac{(\sin^2\theta+b^2\cos^2\theta)^2}{b^2}\sum_{l=1}^\infty\sum_{m=-l,}^l\tilde{\psi}_l^mY_l^m \label{4.11}
\end{equation}

\end{lemma}

\begin{theorem}\thmlab{thm4.1} For any div-free vector field ${u}\in H^k(\mathbb{E}^2)$ and $k\geq0$,

\begin{equation}
    \|{u}-\prod_{\mathrm{null}(\mathcal{L})}{u}\|_{H^k(\mathbb{E}^2)}\leq\|\mathcal{L}[{u}]\|_{H^{k+2}(\mathbb{E}^2)}
\end{equation}

\end{theorem}

\noindent\textbf{Proof of Theorem 4.1.} Consider the stream function $\Psi$, so that $u=\nabla^{\perp}\Psi$, combining with \eqref{4.7} and \eqref{4.11}, we obtain

$$\|{u}-\prod_{\mathrm{null}(\mathcal{L})}{u}\|_{H^k(\mathbb{E}^2)}\sim\sqrt{\sum_{l,m,}\Lambda_{l,m}^{k+1}|\tilde{\psi}_l^m|^2}$$
combining \eqref{4.7} and \eqref{sphericalharmonicrepresentationofL}, we obtain that 

$$\|\mathcal{L}[{u}]\|_{H^{k+2}(\mathbb{E}^2)}=\sqrt{\sum_{\begin{array}{c}l,m,\\m\neq0\end{array}}\Lambda_{l,m}^{k+3}|\frac{m\tilde{\psi}_l^m}{\Lambda_{l,m}}|^2}$$
with observation

$$\Lambda_{l,m}^{k+1}\leq\Lambda_{l,m}^{k+3}\left|\frac{m}{\Lambda_{l,m}}\right|^2$$
Thus one have

$$\|{u}-\prod_{\mathrm{null}(\mathcal{L})}{u}\|_{H^k(\mathbb{S}^2)}\leq\|\mathcal{L}[{u}]\|_{H^{k+2}(\mathbb{E}^2)}$$

\qed

\section{Uniform estimates independent of $\omega$}

In this section, we use energy methods to prove local-in-time existence of classical solutions for the incompressible Euler equation in dependent of $\omega$. Recall the equation \eqref{partialtubyLopertaor}.

\begin{equation}
    \partial_tu+\nabla^\perp\Delta^{-1}\mathrm{curl}(\nabla_uu)=\omega\cdot\mathcal{L}[u] \label{5.1}
\end{equation}
where the operator $\mathcal{L}$ as in \eqref{Loperator} defined by

\begin{equation}
    \mathcal{L}[u]=\nabla^\perp\Delta^{-1}\mathrm{curl}(\frac{-2\sin\theta}{\sqrt{\sin^2\theta+b^2\cos^2\theta}}u^\perp) \label{5.2}
\end{equation}

\noindent The standard energy method (e.g.\cite[Ch.3]{MajdaBertozzi2002}) can by employed to prove. We follow the same process as Cheng\cite{Cheng2013} did. Basically, by the Sobolev embedding theorem on Riemannian manifolds \cite[Ch.2 §10]{Aubin1982} 
\begin{equation}
    |f|_{W^{1,\infty}(\mathbb{E}^2)}\leq C\|f\|_{H^3(\mathbb{E}^2)} \label{5.3}
\end{equation}
To prepare for the proofs in $H^k$ norms, we first show that the
physical energy, $\|{u}\|_{L^2}$ is actually conserved in time,

\begin{equation}
    \|{u}(t,\cdot)\|_{L^2}=\|{u}(0,\cdot)\|_{L^2} \label{5.4}
\end{equation}
Taking the $L^2$ inner product of $u$ with both sides of \eqref{5.1}, we obtain

\begin{equation}
    \langle u,\partial_{t}u\rangle=-\langle u,\nabla^{-1}\Delta^{-1}\mathrm{curl}(\nabla_uu)\rangle+w\langle u,\mathcal{L}[u]\rangle \label{5.5}
\end{equation}
The LHS equals $\frac{1}{2}\partial_t\|u\|_{L^2}^2$, we want to show the RHS is zero. For the second term, 
\begin{equation}
    \begin{aligned}\langle u,\mathcal{L}[u]\rangle&=\langle u,\nabla^{\perp}\Delta^{-1}\mathrm{curl} (\frac{-2\sin\theta}{\sqrt{\sin^{2}\theta+b^{2}\cos^{2}\theta}}u^{\perp})\rangle\\&=\langle\nabla^{\perp}\Delta^{-1}\mathrm{curl} u,\frac{-2\sin\theta}{\sqrt{\sin^{2}\theta+b^{2}\cos^{2}\theta}}u^{\perp}\rangle\\&=\langle u_{\mathrm{inc}},\frac{-2\sin\theta}{\sqrt{\sin^{2}\theta+b^{2}\cos^{2}\theta}}u^{\perp}\rangle=0\end{aligned} \label{5.6}
\end{equation}
By a similar procedure, one can show that $\langle u,\nabla^{\perp}\Delta^{-1}\mathrm{curl}(\nabla_uu)\rangle=\langle u,\nabla_uu\rangle$, to show this part vanishes we shall use \lemref{5.1} which has been proved by \cite{Cheng2013}.

\begin{lemma}\cite{Cheng2013}\lemlab{5.1} Consider vector fields u,v and scalar field fall in $C^1(\mathbb{E}^2)$, and u,v both tangent to $\mathbb{E}^2$.Then
$$\mathrm{div}_{\mathbb{E}^2}{u}=0 \quad\text{implies }\quad\langle f,\nabla_{u}f\rangle=\langle{v},\nabla_{u}{v}\rangle=0$$
\end{lemma}
Combine this lemma and \eqref{5.5}, \eqref{5.6}, we can prove \eqref{5.4}. To make generalization of the above procedure for $\|f{u}\|_{H^k}$, one take spatial derivative $D^{\alpha}$ of \eqref{5.1} up to order $k$, i.e., $|\alpha|\leq k$, and inner-product it with $D^{\alpha}u$, 

\begin{equation}
    \langle D^{\alpha} u,D^{\alpha}\partial_{t}u\rangle=-\langle D^{\alpha}u,D^{\alpha}\nabla^{-1}\Delta^{-1}\mathrm{curl}(\nabla_uu)\rangle+w\langle D^{\alpha}u,D^{\alpha}\mathcal{L}[u]\rangle \label{5.8}
\end{equation}

Now the LHS becomes $\frac{1}{2}\partial_t\|D^\alpha{u}\|_{L^2}^2$ and we need to make sure every term on the RHS is bounded by $\|u\|_{H^{k}}$. We first need to show $\langle D^\alpha{u},D^\alpha\mathcal{L}[{u}]\rangle=\langle D^\alpha{u},\mathcal{L}[D^\alpha{u}]\rangle$. While it is not always be true, only a few set of differential–integral operators on $\mathbb{E}^{2}$ commutes with each other. Let us intruduce the following lemma. Although the conclusion coincides with the spherical case, the proof requires certain modifications due to the geometric differences of the ellipsoidal surface.

\begin{lemma}\cite{Cheng2013} \lemlab{lem5.2} Given integer $j\geq0$.For sufficiently smooth vector field ${u}$ on $\mathbb{E}^2$ with $\mathrm{div}~ u =0$. Then we have 

\begin{equation}\langle\Delta^j{u},\Delta^j\mathcal{L}[{u}]\rangle=0 \quad\text{and}\quad\langle\Delta^j\mathrm{curl}{u},\Delta^j\mathrm{curl}\mathcal{L}[{u}]\rangle=0   
\end{equation}

\end{lemma}

\noindent\textbf{Proof of Lemma 5.2.} It is sufficient to prove  
\begin{equation}
    \langle\Delta^{k}u,\mathcal{L}[u]\rangle=0 \label{5.9}
\end{equation}
due to the symmetric property of $\Delta$, the commutativity of $\mathrm{curl}$, $\langle\nabla^\perp f,{u}\rangle=-\langle f\mathrm{,curl~}{u}\rangle$ and \eqref{hodgeproperty}. We first show $\Delta$ and $\mathcal{L}$ share ellipsoidal harmonics as eigenfunctions and then they commute. 
We have:\\

$\mathcal{L}[u]=\nabla^{\perp}\Delta^{-1}\mathrm{curl}(\frac{-2\sin\theta}{\sqrt{\sin^{2}\theta+b^{2}\cos^{2}\theta}}u^{\perp}) \quad \mathrm{curl}(\frac{-2\sin\theta}{\sqrt{\sin^{2}\theta+b^{2}\cos^{2}\theta}}u^{\perp})=(\nabla\frac{-2\sin\theta}{\sqrt{\sin^{2}\theta+b^{2}\cos^{2}\theta}})u$\\
Thus
$$\mathcal{L}[\nabla^{\perp}\Psi]=\nabla^{\perp}\Delta^{-1}(\nabla\frac{-2\sin\theta}{\sqrt{\sin^{2}\theta+b^{2}\cos^{2}\theta}}\cdot\nabla^{\perp}\Psi)=\nabla^{\perp}\Delta^{-1}[\frac{b^{2}}{(\sin^{2}\theta+b^{2}\cos^{2}\theta)^{2}}\partial_\varphi\Psi]$$
By setting $\Psi=Y_{l}^{m}$ where $l\geq1$, we have
$$\mathcal{L}[\nabla^{\perp}\Psi]=\nabla^{\perp}\sum_{m\neq0}\frac{-im}{\Lambda_{l,m}}\Psi_{l}^{m}Y_{l}^{m} \quad\text{where}\quad\Psi=\widehat{\Psi}_{l}^{m}Y_{l}^{m}$$
So we have $\mathcal{L}[\nabla^{\perp}Y_{l}^{m}]=\frac{-im}{\Lambda_{l,m}}\nabla^{\perp}Y_{l}^{m}$, then 
$$\Delta^{k}\mathcal{L}[\nabla^{\perp}Y_{l}^{m}]=\Delta^{k}\cdot\frac{-im}{\Lambda_{e,m}}\nabla^{\perp}Y_{l}^{m}=\frac{-im}{\Lambda_{l,m}}\Delta^{k}\nabla^{\perp}Y_{l}^{m}=\mathcal{L}[\Delta^{k}\nabla^{\perp}Y_{l}^{m}]$$
Therefore, for any incompressible flow
$$u=\nabla^{\perp}\Psi=\nabla^{\perp}\sum_{l=1}^{\infty}\sum_{m=-l}^{l}\frac{(\sin^{2}\theta+b^{2}\cos\theta)^{2}}{b^{2}}\tilde{\psi}_{l}^{m}Y_{l}^{m}$$ 
we have $\Delta^{k}\mathcal{L}[u]=\mathcal{L}[\Delta^{k}u]$. Combining above properties together we obtain that 
$$\begin{aligned}\langle\Delta^ku,\mathcal{L}[u]\rangle&=-\langle \mathcal{L}[\Delta^ku],u\rangle=-\langle\Delta^k\mathcal{L}[u],u\rangle\\&=-\langle \mathcal{L}[u],\Delta^ku\rangle\end{aligned}$$
which leads to conclusion \eqref{5.9}
\qed\\

\noindent This lemma suggests that we replace $\Delta^k$ with $\Delta^j$ for $k=2j$~(respectively. $\Delta^j\mathrm{curl}$ for $k=2j+1$) in \eqref{5.8}. This is enough to meet the purpose of $H^k$ estimates, because by the definition of $H^k$ norm in \eqref{4.7}, we have 

\begin{equation}
\text{for div }{u}=0,\quad\|{u}\|_{H^{2j}}^2\approx\|\Delta^j{u}\|_{L^2}^2 \quad\text{and}\quad\|{u}\|_{H^{2j+1}}^2\approx\|\Delta^j\operatorname{curl}~{u}\|_{L^2}^2\label{5.12}
\end{equation}
Note that the $\omega$ term in \eqref{5.8} equals to zero, so we only need to estimate nonlinear term, which contains derivatives up to order $k+1$. By \lemref{5.1}, we obtain
\begin{equation}
\quad\langle D^{\alpha}u,D^{\alpha}\nabla^{\perp}\Delta^{-1}\mathrm{curl}\left(\nabla_uu\right)\rangle =\langle D^{\alpha}u,(D^{\alpha}\nabla^{\perp}\Delta^{-1}\mathrm{curl}\left(\nabla_uu\right)-\nabla uD^{\alpha}u)\rangle \label{5.13}
\end{equation}
where $D^\alpha=\Delta^j\mathrm{curl}$ for $k=2j+1$~(resp. $D^\alpha=\Delta^j\mathrm{~for~}k=2j)$.
Note that we have such property $\Delta$ commutes with each one of $\mathrm{div}, \mathrm{curl}, \nabla, \nabla^\perp$. So the commutator term on the RHS of \eqref{5.13} becomes\\

$\begin{aligned}D^{\alpha}\nabla^{\perp}\Delta^{-1}\mathrm{curl}(\nabla_u u)-\nabla uD^{\alpha}u&=\Delta^{j}\mathrm{curl}\nabla^{\perp}(\Delta^{-1}\mathrm{curl}(\nabla_uu))-\nabla u(\Delta^{j}\mathrm{curl}u)\\&=\Delta^{j}\mathrm{curl}(\nabla_uu)-\nabla u(\Delta^{j}\mathrm{curl}u) \quad \text{For $k=2j+1$}
\end{aligned}$\\

resp. $\nabla^{\perp}\Delta^{j+1}\operatorname{curl}(\nabla_{u}u)-\nabla u(\Delta^{j}u)\quad\mathrm{For}\quad k=2j$ \\

\begin{lemma}\cite{Cheng2013}\lemlab{lem5.3} For integer $k=2j+1\geq3$(resp. $k=2j\geq4$) and incompressible vector field $u$ with sufficient regularity, 
$$\|\Delta^j\operatorname{curl}\left(\nabla_{u}{u}\right)-\nabla_{u}(\Delta^j\operatorname{curl}{u})\|_{L^2}\leq C_k\|{u}\|_{H_k}^2$$
resp. $\|\nabla^\perp\Delta^{j-1}\operatorname{curl}\left(\nabla_{u}{u}\right)-\nabla_{u}(\Delta^j{u})\|_{L^2}\leq C_k\|{u}\|_{H_k}^2$

\end{lemma}

\noindent\textbf{Proof of Lemma 5.3.} The key point is that, although the highest derivatives in
each term of the commutator are of $(k+1)$th order, they are
canceled out by the subtraction in the commutator. For the complete calculation, please refer to the \appref{app}. We start from the form of $\mathrm{curl~}(\nabla_{u}{u})-\nabla_{u}(\mathrm{curl~}{u})$, we have 

$$\mathrm{curl~}(\nabla_{{u}}{u})-\nabla_{{u}}(\mathrm{curl~}{u})=\sum_{a,b,c,d=0}^1g_{abcd}(\partial_\theta^au_\theta,\partial_\phi^bu_\theta,\partial_\theta^cu_\phi,\partial_\phi^du_\phi)=\sum_{a,b,c,d=0}^1B(\nabla_3^{a+c}u,\nabla_3^{b+d}u)$$
where $g_{abcd}$ denotes some generic multi-linear function, homogeneous with degree 2 with smooth coefficients. $B(\cdot,\cdot)$ denotes some generic bilinear function with smooth coefficients. Also, $\nabla_3^{a+c}$ denotes $\partial_x,\partial_y,\partial_z$ derivatives and their combinations up to order $a+c$ and these derivatives are taken on the Cartesian components of $u$. Notice that, although 2nd derivatives appear on the LHS, they are canceled out on the RHS.

\noindent Apply $\Delta^{j}$ to above equation, we obtain

\begin{equation}
    \Delta^j\operatorname{curl}\left(\nabla_{u}{u}\right)-\Delta^j\nabla_{u}(\operatorname{curl}{u})=\sum_{\substack{a,b,c,d=0 \\ a+b+c+d\le 2j+2}}^{2j+1}B(\nabla_3^{a+c}u,\nabla_3^{b+d}u) \label{5.14}
\end{equation}
Note that the highest order of derivatives on the LHS is $(2j+2)=(k+1)$, but they are canceled out on the RHS. And $a+b+c+d\le 2j+2$ implies that $\min\{a+c,b+d\}\leq j+1$. We set $a+c=\alpha$, $b+d=\beta$. Thus, the $L^2(\mathbb{E}^2)$ norm of every term on the RHS is bounded by a constant times

$$\sum_{\alpha=0}^{j+1}\sum_{\beta=0}^{2j+1}|\nabla_3^\alpha{u}|_{L^\infty(\mathbb{E}^2)}\|\nabla_3^\beta {u}\|_{L^2(\mathbb{E}^2)}$$

Due to the normally constant extension we did, i.e., the normal component of u and the normal derivatives of u are zero, we could replace  $\nabla_3^\alpha$ or $\nabla_3^\beta$ with combinations of tangential derivatives on $\mathbb{E}^2$. By the standard Sobolev embedding on a compact smooth 2D manifold such that $\mathbb{E}^2$. We could find an upper bound for the above quantity.

$$C|{u}|_{W^{j+1,\infty}(\mathbb{E}^2)}\|{u}\|_{H^{2j+1}(\mathbb{E}^2)}$$

By Sobolev inequality \eqref{5.3}, we have $|{u}|_{W^{j+1,\infty}(\mathbb{S}^2)}$ term of this quantity is bounded by the $\|{u}\|_{H^{2j+1}(\mathbb{S}^2)}$ term, as long as $j\geq1$. Thus we could establish an upper bound for the $L^2(\mathbb{E}^2)$ norm of \eqref{5.14}.

\begin{equation}
    \|\Delta^j\operatorname{curl}(\nabla_{u}{u})-\Delta^j\nabla_{u}(\operatorname{curl}{u})\|_{L^2(\mathbb{E}^2)}\leq C\|{u}\|_{H^k}^2 \label{5.15}
\end{equation}

The same type of calculation works for estimating the $L^2$ norm of commutator $\Delta^j\nabla_{u}(\operatorname{curl}{u})-\nabla_{u}(\Delta^j\operatorname{curl}{u})$. Thus the same type of bound
in \eqref{5.15} also works for this term, so we have:

$$\|\Delta^j\nabla_{u}(\operatorname{curl}{u})-\nabla_{u}(\Delta^j\operatorname{curl}{u})\|_{L^2(\mathbb{E}^2)}\leq C\|{u}\|_{H^k}^2$$

Adding up two estimates and applying the triangle inequality, we prove the first part of the conclusion. The second part
for $k=2j\geq4$ can be proved in similar way. 

\qed

\noindent After we get the bound of commutator, we could obtain the theorem for $H^k$ estimates of the solution independent of  $\omega$.

\begin{theorem}\thmlab{thm5.1}
Consider the incompressible Euler equations \eqref{5.1}, \eqref{5.2} on a rotating ellipsoid $\mathbb{E}^2$ with div-free initial data $u_0.$ Given any integer $k\geq3$, assume $u_0\in H^k(E^2).$ Then, there exists universal constants $C_0,T_0$ independent of $\omega$ so that

$$\|{u}(t,\cdot)\|_{H^k}\leq C_0\|{u}_0\|_{H^k}\quad\text{for any }t\in\left[0,\frac{T_0}{\|{u}_0\|_{H^k}}\right]$$   
\end{theorem}

\noindent\textbf{Proof of Theorem 5.1.} The existence and uniqueness of $H^k$ solution is well established for general hyperbolic PDE systems that are symmetrizable, e.g. \cite{Taylor1997}. Assume the maximum life span of such solution is $T_\omega$. Our goal is to show that $T_{\omega}\geq T_{0}/\|{u}_{0}\|_{H^{k}}$.\\

Set $D^\alpha=\Delta^j\operatorname{curl}$ for $k=2j+1$ (resp. $D^\alpha=\Delta^j\operatorname{for}k=2j$) in \eqref{5.8}. Apply \lemref{lem5.2} to cancel out the $\omega$ term and apply \eqref{5.13} with \lemref{lem5.3} and Cauchy-Schwartz inequality to estimate tri-linear product,

$$\frac{1}{2}\partial_{t}\|D^{\alpha}{u}\|_{L^{2}}^{2}\leq C\|D^{\alpha}{u}\|_{L^{2}}\|{u}\|_{H^{k}}^{2}$$
By equivalence of $H^K$ norms as in \eqref{5.12}, the above estimate becomes to 

$$\frac{1}{2}\partial_t\|D^\alpha{u}\|_{L^2}^2\leq\mathbb{C}\|D^\alpha{u}\|_{L^2}^3$$

By the uniqueness of classical solutions,${u}_0\equiv0\Longleftrightarrow{u}(t,\cdot)\equiv0$; thus, we only deal with solution $u$ with $\|D^\alpha{u}\|_L^2\approx\|{u}\|_{H^k}\neq0.$ Simplify the above estimate by dividing both sides with $\left\|D^\alpha{u}\right\|_{L^2}^3$

$$\frac{\partial_t\|D^\alpha{u}\|_{L^2}}{\|D^\alpha{u}\|_{L^2}^2}\leq C\Longrightarrow-\partial_t\left(\|D^\alpha{u}\|_{L^2}^{-1}\right)\leq C\Longrightarrow\|D^\alpha{u}\|_{L^2}\leq\frac{1}{\|D^\alpha{u}_0\|_{L^2}-Ct}$$

By the equivalence \eqref{5.12}, we can find suitable values for $C_0$, $T_0$ as used the conclusion of the theorem. Note that all constants in this proof are independent of $\omega$ and only depent on $k$.

\section{Proof of the main theorem}

The main \thmref{thm1.1} fits into the framework of \lemref{lemma1.1}. And we have obtained most ingredients in this framework. We have defined and talked about properties of operator $\mathcal{L}$ in \defref{def2.1}, \eqref{identity} and \lemref{lem4.1}; and defined operator $\prod_{\mathrm{null}}$ in \lemref{lem3.2}. The key estimate \eqref{1.7} has been proved in \thmref{thm4.1} and $\omega$-independent $H^k$ estimates for the solution in \thmref{thm5.1}. Now the only thing left are regularity assumptions \eqref{1.5} and \eqref{1.6}. We shall prove them and establish the estimates for $M_0$ in \eqref{1.3}.

First by fitting \ref{5.1} into the framework \lemref{lemma1.1}, we set $\mathcal{L}$ as defined in \eqref{5.2} and $f:=\nabla^\perp\Delta^{-1}\operatorname{curl}\left(\nabla_{u}{u}\right)$.
Given initial data $u_0 \in H^k(\mathbb{S}^2)$, by \thmref{thm5.1}, we have 

\begin{equation}
    \|(t,\cdot)\|_{H^k}\leq C_0\|{u}_0\|_{H^k}\quad\text{for any }t\in\left[0,\frac{T_0}{\|{u}_0\|_{H^k}}\right] \label{6.1}
\end{equation}
Then, by the proof of \lemref{lem5.3}, we can show that 

\begin{equation}
    \|f\|_{H^{k-1}}\leq C\|{u}\|_{H^k}^2\quad\mathrm{and}\quad\|\mathcal{L}[{u}]\|_{H^{k+1}}\leq C\|{u}\|_{H^k} \label{6.2}
\end{equation}
Set Hilbert spaces: $X_1:=H^{k-3}(\mathbb{E}^2)$, $X_2:=H^{k-1}(\mathbb{E}^2)$. Assumption \eqref{1.7} is verified by \thmref{thm4.1}, assumption \eqref{1.5} is verified by \eqref{6.1}, \eqref{6.2}. Note that the time-continuity part of \eqref{1.5} is due to the calculation,

$$|\|f{u}(t+\tau,\cdot)-{u}(t,\cdot)\|_{H^k}|\leq|\|{u}(t+\tau,\cdot)\|_{H^k}-\|{u}(t,\cdot)\|_{H^k}|=\left|\int_t^{t+\tau}\partial_t\|{u}\|_{H^k}\right|$$
and the $\partial_t\|{u}\|_{H^k}$ term on the RHS can be shown to be uniformly bounded following the proof of \thmref{thm5.1}. Then, one lets $\tau\to\tilde{0}$ to show the time-continuity of $u$. Likewise, this can be done for $\mathcal{L}[{u}]$ and$f$.

As for the commutativity of $\mathcal{L}$ and the integral operator in time in \ref{1.6}, we use the time continuity of $u$ and $\mathcal{L}[f{u}]$ in terms of $H^k$ norm to rewrite the time integrals as limits of Riemann sums

$$\int_0^T\mathcal{L}[{u}]dt=\lim_{N\to0}\sum_i\mathcal{L}[{u}(t_i,\cdot)]\delta t$$
$$\mathcal{L}\left[\int_0^T{u}\:dt\right]=\mathcal{L}\left[\lim_{N\to0}\sum_i{u}(t_i,\cdot)\delta t\right]$$
where $\delta t=T/N$ and $0=t_0<t_1<t_2<\cdots<t_N=T$ form an equi-partition of $[0,T]$. In the RHS of the second equality, $\mathcal{L}$ and lim commute because $\mathcal{L}$ is a continuous mapping on $H^k$. Thus, the LHS of the two above equalities are equal.

Finally, by \eqref{6.1}, \eqref{6.2}, the constant $M$ in \lemref{lemma1.1} is bounded by $\|{u}\|_{H^k}\leq C_0\|{u}_0\|_{H^k}$ and the constant $M^\prime$ is bounded by $\left\|f\right\|_{H^{k-1}}\leq C\left\|{u}\right\|_{H^k}^2\leq C^{\prime}\left\|{u}_0\right\|_{H^k}^2$. These two bounds validate the use of constant $M_0$ in \ref{1.3}. 

The proof of \thmref{thm1.1} is complete.

\section*{Appendix}\label{app}

\noindent\textbf{Proof of \lemref{lem5.3}} Denote $\mathbb{E}^2$ as the ellipsoid.

$$x=\rho\cos\theta\cos\phi,\quad y=\rho\cos\theta\sin\phi,\quad z=b\rho\sin\theta.$$
Set

$$\begin{aligned}&\mathbf{e}_\rho=(\cos\theta\cos\phi,\cos\theta\sin\phi,b\sin\theta),\\&\mathbf{e}_\phi=(-\sin\phi,\cos\phi,0),\\&\mathbf{e}_\theta=(-\sin\theta\cos\phi,-\sin\theta\sin\phi,b\cos\theta).\end{aligned}$$
This is a non-orthogonal non-normalized right-handed frame for $T\mathbb{R}^3$ at any point
except the origin. Also,
$$\begin{aligned}&\mathbf{e}_x=\cos\theta\cos\phi\mathbf{e}_\rho-\sin\phi\mathbf{e}_\phi-\sin\theta\cos\phi\mathbf{e}_\theta,\\&\mathbf{e}_y=\cos\theta\sin\phi\mathbf{e}_\rho+\cos\phi\mathbf{e}_\phi-\sin\theta\sin\phi\mathbf{e}_\theta,\\&\mathbf{e}_z=\frac1b\sin\theta\mathbf{e}_\rho+\frac1b\cos\theta\mathbf{e}_\theta.\end{aligned}$$
For

$$\begin{aligned}
    \mathbf{v}&=v_x\mathbf{e}_x+v_y\mathbf{e}_y+v_z\mathbf{e}_z=v_\rho\mathbf{e}_\rho/\|\mathbf{e}_\rho\|+v_\phi\mathbf{e}_\phi/\|\mathbf{e}_\phi\|+v_\theta\mathbf{e}_\theta/\|\mathbf{e}_\theta\|\\&=v_\rho\frac{1}{\sqrt{\cos^2\theta+b^2\sin^2\theta}}\partial_\rho+v_\phi\frac{1}{\cos\theta}\partial_\phi+v_\theta\frac{1}{\sqrt{\sin^2\theta+b^2\cos^2\theta}}\partial_\theta
\end{aligned}$$
One has\\

$\begin{aligned}&v_x=(v_\rho\frac{\mathbf{e}_\rho}{\|\mathbf{e}_\rho\|}+v_\phi\frac{\mathbf{e}_\phi}{\|\mathbf{e}_\phi\|}+v_\theta\frac{\mathbf{e}_\theta}{\|\mathbf{e}_\theta\|})\cdot\mathbf{e}_x=v_\rho\cos\theta\cos\phi/\|\mathbf{e}_\rho\|-v_\phi\sin\phi/\|\mathbf{e}_\phi\|-v_\theta\sin\theta\cos\phi/\|\mathbf{e}_\theta\|,\\&v_y=(v_\rho\frac{\mathbf{e}_\rho}{\|\mathbf{e}_\rho\|}+v_\phi\frac{\mathbf{e}_\phi}{\|\mathbf{e}_\phi\|}+v_\theta\frac{\mathbf{e}_\theta}{\|\mathbf{e}_\theta\|})\cdot\mathbf{e}_y=v_\rho\cos\theta\sin\phi/\|\mathbf{e}_\rho\|+v_\phi\cos\phi/\|\mathbf{e}_\phi\|-v_\theta\sin\theta\sin\phi/\|\mathbf{e}_\theta\|,\\&v_z=(v_\rho\frac{\mathbf{e}_\rho}{\|\mathbf{e}_\rho\|}+v_\phi\frac{\mathbf{e}_\phi}{\|\mathbf{e}_\phi\|}+v_\theta\frac{\mathbf{e}_\theta}{\|\mathbf{e}_\theta\|})\cdot\mathbf{e}_z=bv_\rho\sin\theta/\|\mathbf{e}_\rho\|+bv_\theta\cos\theta/\|\mathbf{e}_\theta\|\end{aligned}$
Now if $\mathbf{v}\in T\mathbb{E}^2$, $v_\rho=0$:

$$\begin{aligned}&v_{x}=-v_\phi\sin\phi/\|\mathbf{e}_\phi\|-v_\theta\sin\theta\cos\phi/\|\mathbf{e}_\theta\|,\\&v_{y}=v_\phi\cos\phi/\|\mathbf{e}_\phi\|-v_\theta\sin\theta\sin\phi/\|\mathbf{e}_\theta\|,\\&v_{z}=bv_\theta\cos\theta/\|\mathbf{e}_\theta\|.\end{aligned}$$

$\begin{aligned}&\nabla_{\mathbf{v}}\mathbf{v}\\&=\underset{T\mathbb{R}^3\to T\mathbb{E}^2}{\operatorname*{\operatorname*{Proj}}}\sum_{j=x,y,z}(\nabla_\mathbf{v}v_i)\mathbf{e}_i\\&=\underset{T\mathbb{R}^3\to T\mathbb{E}^2}{\operatorname*{\operatorname*{Proj}}}[\nabla_{\mathbf{v}}(-v_\phi\sin\phi/\|\mathbf{e}_\phi\|-v_\theta\sin\theta\cos\phi/\|\mathbf{e}_\theta\|)(\cos\theta\cos\phi\mathbf{e}_\rho-\sin\phi\mathbf{e}_\phi-\sin\theta\cos\phi)\\&+\nabla_{\mathbf{v}}(v_\phi\cos\phi/\|\mathbf{e}_\phi\|-v_\theta\sin\theta\sin\phi/\|\mathbf{e}_\theta\|)(\cos\theta\sin\phi\mathbf{e}_\rho+\cos\phi\mathbf{e}_\phi-\sin\theta\sin\phi\mathbf{e}_\theta)\\&+\nabla_\mathbf{v}(bv_\theta\cos\theta/\|\mathbf{e}_\theta\|)(\frac{1}{b}\sin\theta\mathbf{e}_\rho+\frac{1}{b}\cos\theta\mathbf{e}_\theta)\\&=\mathbf{e}_\phi[(v_\phi\frac{1}{\cos\theta}\frac{\partial}{\partial\phi}+v_\theta\frac{1}{\sqrt{\sin^2\theta+b^2\cos^2\theta}}\frac{\partial}{\partial\theta})(v_\phi\sin\phi/\|\mathbf{e}_\phi\|+v_\theta\sin\theta\cos\phi/\|\mathbf{e}_\theta\|)\sin\phi\\&+(v_\phi\frac{1}{\cos\theta}\frac{\partial}{\partial\phi}+v_\theta\frac{1}{\sqrt{\sin^2\theta+b^2\cos^2\theta}}\frac{\partial}{\partial\theta})(v_\phi\cos\phi/\|\mathbf{e}_\phi\|-v_\theta\sin\theta\sin\phi/\|\mathbf{e}_\theta\|)\cos\phi]\\&+\mathbf{e}_\theta[(v_\phi\frac{1}{\cos\theta}\frac{\partial}{\partial\phi}+v_\theta\frac{1}{\sqrt{\sin^2\theta+b^2\cos^2\theta}}\frac{\partial}{\partial\theta})(v_\phi\sin\phi/\|\mathbf{e}_\phi\|+v_\theta\sin\theta\cos\phi/\|\mathbf{e}_\theta\|)\sin\theta\cos\phi\\&-(v_\phi\frac{1}{\cos\theta}\frac{\partial}{\partial\phi}+v_\theta\frac{1}{\sqrt{\sin^2\theta+b^2\cos^2\theta}}\frac{\partial}{\partial\theta})(v_\phi\cos\phi/\|\mathbf{e}_\phi\|-v_\theta\sin\theta\sin\phi/\|\mathbf{e}_\theta\|)\sin\theta\sin\phi\\&+(v_\phi\frac{1}{\cos\theta}\frac{\partial}{\partial\phi}+v_\theta\frac{1}{\sqrt{\sin^2\theta+b^2\cos^2\theta}}\frac{\partial}{\partial\theta})(bv_\theta\cos\theta/\|\mathbf{e}_\theta\|)\frac{1}{b}\cos\theta]
\end{aligned}$\\

For ${u}\in T\mathbb{E}^2$:

$${u}=u_\phi\mathbf{e}_\phi/\|\mathbf{e}_\phi\|+u_\theta\mathbf{e}_\theta/\|\mathbf{e}_\theta\|=u_\phi\frac{1}{\cos\theta}\frac{\partial}{\partial\phi}+u_\theta\frac{1}{\sqrt{\sin^2\theta+b^2\cos^2\theta}}\frac{\partial}{\partial\theta}$$

$$\mathrm{curl~}{u}=\mathrm{div~}{u}^\perp=-\frac{1}{\cos\theta}\partial_\phi u_\theta+\frac{1}{\cos\theta\sqrt{\sin^2\theta+b^2\cos^2\theta}}\partial_\theta(\cos\theta u_\phi)$$
Hence\\

$\begin{aligned}&\mathrm{curl~}(\nabla_{u}\mathbf{u})\\&=-\frac{1}{\cos\theta}\partial_\phi[\|\mathbf{e}_\theta\|(u_\phi\frac{1}{\cos\theta}\frac{\partial}{\partial\phi}+u_\theta\frac{1}{\sqrt{\sin^2\theta+b^2\cos^2\theta}}\frac{\partial}{\partial\theta})(u_\phi\sin\phi/\|\mathbf{e}_\phi\|+u_\theta\sin\theta\cos\phi/\|\mathbf{e}_\theta\|)\sin\theta\cos\phi\\&-(u_\phi\frac{1}{\cos\theta}\frac{\partial}{\partial\phi}+u_\theta\frac{1}{\sqrt{\sin^2\theta+b^2\cos^2\theta}}\frac{\partial}{\partial\theta})(u_\phi\cos\phi/\|\mathbf{e}_\phi\|-u_\theta\sin\theta\sin\phi/\|\mathbf{e}_\theta\|)\sin\theta\sin\phi\\&+(u_\phi\frac{1}{\cos\theta}\frac{\partial}{\partial\phi}+v_\theta\frac{1}{\sqrt{\sin^2\theta+b^2\cos^2\theta}}\frac{\partial}{\partial\theta})(bu_\theta\cos\theta/\|\mathbf{e}_\theta\|)\frac{1}{b}\cos\theta]\\&+\frac{1}{\cos\theta\sqrt{\sin^2\theta+b^2\cos^2\theta}}\\&\cdot\partial_\theta(\cos\theta\|\mathbf{e}_\phi\|[(u_\phi\frac{1}{\cos\theta}\frac{\partial}{\partial\phi}+u_\theta\frac{1}{\sqrt{\sin^2\theta+b^2\cos^2\theta}}\frac{\partial}{\partial\theta})(u_\phi\sin\phi/\|\mathbf{e}_\phi\|+u_\theta\sin\theta\cos\phi/\|\mathbf{e}_\theta\|)\sin\phi\\&+(u_\phi\frac{1}{\cos\theta}\frac{\partial}{\partial\phi}+u_\theta\frac{1}{\sqrt{\sin^2\theta+b^2\cos^2\theta}}\frac{\partial}{\partial\theta})(u_\phi\cos\phi/\|\mathbf{e}_\phi\|-u_\theta\sin\theta\sin\phi/\|\mathbf{e}_\theta\|)\cos\phi]).\end{aligned}$\\
Also\\

$\begin{aligned}&\nabla_{\mathbf{u}}(\mathrm{curl~}{u})\\&=\nabla_\mathbf{u}(-\frac{1}{\cos\theta}\partial_\phi u_\theta+\frac{1}{\cos\theta\sqrt{\sin^2\theta+b^2\cos^2\theta}}\partial_\theta(\cos\theta u_\phi))\\&=(u_\phi\frac{1}{\cos\theta}\frac{\partial}{\partial\phi}+u_\theta\frac{1}{\sqrt{\sin^2\theta+b^2\cos^2\theta}}\frac{\partial}{\partial\theta})(-\frac{1}{\cos\theta}\partial_\phi u_\theta+\frac{1}{\cos\theta\sqrt{\sin^2\theta+b^2\cos^2\theta}}\partial_\theta(\cos\theta u_\phi)).\end{aligned}$

Combining the results above, with computation we can deduce that the commutator satisfies

$$\mathrm{curl~}(\nabla_{u}{u})-\nabla_{u}(\mathrm{curl~}{u})=\sum_{a,b,c,d=0}^1g_{abcd}(\partial_\theta^au_\theta,\partial_\phi^bu_\theta,\partial_\theta^cu_\phi,\partial_\phi^du_\phi)$$
where $g_{abcd}$ denotes some generic multi-linear function, homogeneous with degree 2,
with smooth coefficients, since the terms with the highest order derivatives cancel out.

\printbibliography

@article{Cheng2013,
  author  = {Cheng, Bin and Mahalov, Alex},
  title   = {Euler equation on a fast rotating sphere—Time-averages and zonal flows},
  journal = {European Journal of Mechanics B/Fluids},
  volume  = {37},
  pages   = {48--58},
  year    = {2013}
}

@article{Xu2024,
  author  = {Xu, Chenghao},
  title   = {The Non-zonal {Rossby--Haurwitz} Solutions of the 2{D} {Euler} Equations on a Rotating Ellipsoid},
  journal = {Journal of Mathematical Fluid Mechanics},
  volume  = {26},
  pages   = {64},
  year    = {2024}
}

@article{Cheng2012_SIAM,
  author    = {Bin Cheng},
  title     = {Singular limits and convergence rates of compressible {Euler} and rotating shallow water equations},
  journal   = {SIAM Journal on Mathematical Analysis},
  volume    = {44},
  number    = {2},
  pages     = {1050--1076},
  year      = {2012},
  doi       = {10.1137/11085147X}
}

@article{ChengMahalov2012_DCDS,
  author    = {Bin Cheng and Alex Mahalov},
  title     = {Time averages of fast oscillatory systems},
  journal   = {Discrete and Continuous Dynamical Systems - Series S},
  year      = {2012},
  note      = {in press}
}

@book{Warner1983,
  author    = {Frank W. Warner},
  title     = {{Foundations of Differentiable Manifolds and Lie Groups}},
  series    = {Graduate Texts in Mathematics},
  volume    = {94},
  publisher = {Springer-Verlag},
  address   = {New York, Berlin},
  year      = {1983},
  note      = {Corrected reprint of the 1971 edition}
}

@book{Taylor1996,
  author    = {Michael E. Taylor},
  title     = {{Partial Differential Equations. I. Basic Theory}},
  series    = {Applied Mathematical Sciences},
  volume    = {115},
  publisher = {Springer-Verlag},
  address   = {New York},
  year      = {1996}
}

@book{MajdaBertozzi2002,
  author    = {Andrew J. Majda and Andrea L. Bertozzi},
  title     = {{Vorticity and Incompressible Flow}},
  series    = {Cambridge Texts in Applied Mathematics},
  volume    = {27},
  publisher = {Cambridge University Press},
  address   = {Cambridge},
  year      = {2002},
  pages     = {xii+545}
}

@book{Aubin1982,
  author    = {Thierry Aubin},
  title     = {{Nonlinear Analysis on Manifolds. Monge--Ampère Equations}},
  series    = {Grundlehren der Mathematischen Wissenschaften (Fundamental Principles of Mathematical Sciences)},
  volume    = {252},
  publisher = {Springer-Verlag},
  address   = {New York},
  year      = {1982},
  pages     = {xii+204}
}

@book{Taylor1997,
  author    = {Michael E. Taylor},
  title     = {{Partial Differential Equations. III. Nonlinear Equations}},
  series    = {Applied Mathematical Sciences},
  volume    = {117},
  publisher = {Springer-Verlag},
  address   = {New York},
  year      = {1997}
}

@article{BerardoWit2022,
  author  = {Berardo, D. and Wit, J.},
  title   = {{On the Effects of Planetary Oblateness on Exoplanet Studies}},
  journal = {The Astrophysical Journal},
  volume  = {935},
  number  = {2},
  pages   = {178},
  year    = {2022},
  doi     = {10.3847/1538-4357/ac7d59}
}

@book{ElkinsTanton2006,
  author    = {Elkins-Tanton, L. T.},
  title     = {{Jupiter and Saturn}},
  publisher = {Infobase Publishing},
  address   = {New York},
  year      = {2006}
}

@article{Haziot2022,
  author  = {Haziot, S. V.},
  title   = {{On the Spherical Geopotential Approximation for Saturn}},
  journal = {Communications on Pure and Applied Analysis},
  volume  = {21},
  number  = {7},
  pages   = {2327},
  year    = {2022},
  doi     = {10.3934/cpaa.2022104}
}

@article{Taylor2016,
  author  = {Michael E. Taylor},
  title   = {Euler Equation on a Rotating Surface},
  journal = {Journal of Functional Analysis},
  volume  = {270},
  number  = {10},
  pages   = {3884--3945},
  year    = {2016},
  doi     = {10.1016/j.jfa.2016.02.014}
}
\end{document}